\documentclass[a4paper, 12pt]{article}
\usepackage[utf8]{inputenc}
\usepackage[font=small,labelfont=bf]{caption}
\usepackage[skip=0.4\baselineskip]{caption}
\usepackage{rotating}
\usepackage{longtable}
\usepackage{amsmath,amsfonts,amssymb}
\usepackage{mathtools}

\usepackage{multirow}
\usepackage{fullpage}
\usepackage{times}
\usepackage{adjustbox}
\usepackage{pgfplots}
\usepackage{tablefootnote}
\usepackage{fancyhdr,graphicx,amsmath,amssymb}
\usepackage[english]{babel}
\usepackage[T1]{fontenc}
\usepackage{placeins}
\usepackage{listings}

\usepackage{setspace}
\usepackage{amssymb}
\usepackage{bbm}
\usepackage{algorithm}
\usepackage{algpseudocode}
\algrenewcommand\algorithmicrequire{\textbf{Input:}}
\algrenewcommand\algorithmicensure{\textbf{Output:}}
\usepackage{accents}

\usepackage{makecell}
\usepackage{amsthm}

\theoremstyle{plain}

\theoremstyle{definition}
 
\theoremstyle{remark}

\algdef{SE}[DOWHILE]{Do}{doWhile}{\algorithmicdo}[1]{\algorithmicwhile\ #1}%

\let\emptyset\varnothing

\algnewcommand{\algorithmicor}{{ or }}
\algnewcommand{\OR}{\algorithmicor}

\usepackage[a4paper,margin=0.9in]{geometry}
\usepackage{eurosym}
\usepackage{amsmath}
\usepackage{natbib}
\usepackage{graphicx}
\usepackage[colorinlistoftodos]{todonotes}
\usepackage[colorlinks=true, allcolors=blue]{hyperref}
\usepackage{authblk}
\usepackage{cleveref}
\usepackage[super]{nth}
\usepackage{enumitem, float, subcaption}
\usepackage[compact]{titlesec}
\usepackage{booktabs,caption}
\usepackage[flushleft]{threeparttable}
\renewcommand{\footnotesize}{\fontsize{8.5pt}{11pt}\selectfont}

\newcommand{\matr}[1]{\mathbf{#1}}

\makeatletter
\newenvironment{smallalign*}
{\par
\footnotesize
\setlength{\abovedisplayskip}{0pt} 
\csname align*\endcsname}
{\csname endalign*\endcsname}
\makeatother

\usepackage{changes}
\usepackage{lipsum}
\definechangesauthor[name={Rick Willemsen}, color=orange]{rick}

\titlespacing*{\section}
{0pt}{0.5em}{0em}  

\titlespacing*{\subsection}
{0pt}{0.5em}{0em}  

\titlespacing*{\subsubsection}
{0pt}{0.5em}{0em}  

\begin{document}
\hypersetup{
    colorlinks = true, 
    citecolor = black,
    linkcolor = black,
    urlcolor = black
}
\title{A Benders decomposition approach for hierarchical $k$-medoids}
\author{Rick S. H. Willemsen}
\affil{Econometric Institute, Erasmus University, 3062 PA Rotterdam, The Netherlands}

\date{\today}

\maketitle

\renewenvironment{abstract}
{\begin{quote}
\noindent \rule{\linewidth}{.5pt}}
{\smallskip\noindent \rule{\linewidth}{.5pt}
\end{quote}
}
\begin{abstract} \\
We study a hierarchical clustering problem formulated as a sum of $k$-medoids objectives. This problem is related to the multi-period $p$-median problem, but instead of having \textit{nested facilities}, we consider a hierarchy of \textit{nested assignments}. We present two Benders formulations from the $p$-median literature and show how the obtained nested $p$-median solutions can be converted into hierarchical $k$-medoids solutions. We validate our methods on various instances. First, we show on small-scale instances that a standard Benders decomposition approach achieves speedups of up to three orders of magnitude over an existing branch-and-price method, while producing hierarchical $k$-medoids solutions that are at most 5\% away from optimality. Second, we present an aggregated Benders heuristic that can handle instances with more than 13,000 observations. \\
\end{abstract}

\section{Introduction}
Cluster analysis deals with partitioning groups of objects based on certain measurements \citep{hennig2015}. Cluster analysis can be divided into hierarchical and nonhierarchical (also known as partitional) clustering \citep{saxena2017}. The focus of this paper is on the hierarchical variant, which is often solved with heuristics. These heuristics are commonly divided into agglomerative and divisive ones and typically employ greedy, local decision rules based on similarity measures and linkage criteria \citep{saxena2017}. 

\paragraph{Background and literature}

Recently, global objectives have been proposed for evaluating clustering hierarchies. Several classes of objective functions have been studied in the literature \citep{dasgupta2016, moseley2017, cohenaddad2019}, along with a number of specific objectives \citep{wang2020, greenberg2021, vichi2022}. These objectives also motivated the study of approximation algorithms \citep{roy2016, charikar2017, cohenaddad2019, wang2020, naumov2021, rajagopalan2021, chami2020}. 

Approaches based on mathematical programming formulations have been studied by \cite{gilpin2013}, \cite{gilpin2017}, and \cite{vichi2022}, while exact solution methods have been proposed by \cite{greenberg2021}, \cite{willemsen2026exact}, \cite{labbe2026mixed}, and \cite{amorosi2026mathematical}. Specifically, \cite{willemsen2026exact} solve an objective defined as a sum of partitional clustering objectives over each level of the hierarchy, generalizing the objective proposed by \cite{vichi2022}, which considers the sum of sum-of-squares clustering objectives. Their exact method is able to solve instances with up to 200 observations to optimality, although the largest instance requires more than 13 hours of computation. They consider the top nine layers of the hierarchy rather than the full hierarchy. 

\paragraph{Relationship between $k$-medoids and $p$-median}
In this paper, we consider the hierarchical $k$-medoids problem, which we define as the sum of $k$-medoids objectives across the levels of the hierarchy, subject to nested cluster assignments. This formulation is a special case of the objective considered by \cite{willemsen2026exact}. Given the close connection between $k$-medoids and $p$-median problems, we propose to leverage techniques from the $p$-median literature to solve the hierarchical $k$-medoids problem.

The $k$-medoids and $p$-median problems are often regarded as equivalent \citep{brusco2019k, kazakovtsev2020application, schubert2021fast, willemsen2025accelerating}. In particular, the $p$-median is a special case of the facility location problem in which there are no costs for opening facilities and exactly $p$ facilities must be selected from the set of existing locations. Under these conditions, the $p$-median problem is equivalent to the $k$-medoids problem, with $p=k$.

However, the equivalence between the $k$-medoids and $p$-median problems does not necessarily extend to a nested setting. For example, consider a nested $p$-median problem arising in a multi-period setting, where the facilities selected in one period are required to remain open in subsequent periods. Although this imposes a nesting constraint on the facilities, it is not equivalent to the hierarchical $k$-medoids problem, which requires the cluster assignments to be nested across levels.

The relationship between nested $p$-median problems and hierarchical $k$-medoids problems can be summarized by the taxonomy in \Cref{tab:taxonomy}. When neither cluster assignments nor center sets are required to be nested, the problem decomposes into independent $k$-medoids or $p$-median problems, which can be solved separately at each level. Requiring nested cluster assignments leads to the hierarchical $k$-medoids problem, whereas enforcing nested center sets results in the nested $p$-median problem. Finally, when both center sets and cluster assignments are nested we call it fully nested $k$-medoids.

\begin{table}[h]
    \centering
    \caption{Taxonomy of nested center sets and cluster assignments.}
    \label{tab:taxonomy}
    \begin{tabular}{lll}
    \toprule
       nested center sets   & \multicolumn{2}{c}{nested cluster assignments} \\  \cline{2-3}
         & not enforced & enforced \\ 
         \midrule
   not enforced      & independent $k$-medoids (or $p$-median) & hierarchical $k$-medoids\\
   enforced      & nested $p$-median & fully nested $k$-medoids \\ \bottomrule
    \end{tabular}
\end{table}

Benders decomposition has been successfully applied to single-level $p$-median and related problems \citep{duranmateluna2023, willemsen2025accelerating, fischetti2016}. Related decomposition methods have also been developed for multi-period facility location problems, where decisions across periods are coupled through temporal constraints  \citep{sonmez2012decomposition, arulselvan2019incremental, ortiz2019exact}, although other mathematical programming based approaches are also used \citep{correia2016multi}. 

Nestedness has also been explicitly incorporated into mathematical programming formulations for related facility location problems. Closely related to our setting, \cite{mcgarvey2022nested} introduce the concept of nested facility locations and apply it to the $p$-median problem. They consider regret-minimizing objectives and develop Lagrangian-relaxation based heuristics for solving the problem. \cite{brandstetter2024nested} study a nested $p$-center problem with objectives minimizing regret as well.

\paragraph{Research gap and contributions}
The branch-and-price method for solving hierarchical $k$-medoids proposed in \cite{willemsen2026exact} has limited scalability. In particular, they solve instances with up to 200 observations to optimality, while the largest instance requires more than 13 hours of computation. Moreover, the reported experiments consider hierarchies with at most nine levels, rather than complete hierarchies. These limitations leave open the question of whether more scalable solution approaches can be developed for the hierarchical $k$-medoids problem.

A second gap concerns the relationship between $k$-medoids and $p$-median methods in the hierarchical setting. Although the two problems are equivalent at a single level, it is not clear to what extent solution methods developed for the $p$-median problem can be adapted to solve hierarchical $k$-medoids problems. To the best of our knowledge, this question has not been systematically investigated in the literature.

Our contributions are as follows. First, we adapt two Benders formulations from the $p$-median literature and convert the obtained solutions into hierarchical $k$-medoids solutions. Second, in our computational experiments we show that a standard Benders decomposition approach can be three orders of magnitude faster than the branch-and-price method for hierarchical $k$-medoids. Third, we show that an aggregated Benders heuristic can solve instances with more than 13,000 observations.

The remainder of this paper is organized as follows. \Cref{section2:problemDefinition} formally introduces the hierarchical $k$-medoids problem. \Cref{section3:formulations} presents the Benders formulations, followed by the solution methods in \Cref{section4:solutionMethods}. Computational results are discussed in \Cref{section5:results}, and \Cref{section6:conclusion} concludes the paper.

\section{Hierarchical $k$-medoids problem}\label{section2:problemDefinition}

Let $O=\{o_1, \dots, o_N\}$ be a set of $N$ objects in $m$-dimensional space. At each level $k$, a partition $P_k=\{C_{k1},\dots, C_{kk}\}$ is needed, which contains $k$ clusters. The clusters satisfy $C_{ki}\neq \emptyset, \forall i\in [k], C_{ki}\cap C_{kj}=\emptyset, \forall i\neq j$ and $\cup_{i=1}^{k} C_{ki} = O$, where we use the convention that $[N]$ is the set $\{1, 2, \dots, N\}$ and $[0, N]$ is the set $\{0, 1, \dots, N\}$. Each cluster has a medoid, a representative object within the cluster that minimizes the sum of distances to all objects in that cluster \citep{kaufman1990}.

Instead of obtaining a full hierarchy consisting of $N$ levels, it is possible to use a parsimonious tree to improve interpretability by selecting the top $K$ levels, where $K\leq N$ \citep{gordon1987}. Note that $|P_k|=k, \forall k\in [K]$. Then, the goal in hierarchical clustering is to obtain $K$ nested partitions $P_1, \dots, P_K$ of the same objects minimizing a given objective function, where the partitions are nested if, for every $k\in [K-1]$, each cluster in $P_{k+1}$ is a subset of some cluster in $P_k$.

In this paper, we consider a special case of the objective in \cite{willemsen2026exact}, namely the sum of the $k$-medoids objectives across all levels in the hierarchy. Let $d_{ij}$ denote the distance between objects $o_i$ and $o_j$, and let $\mu_k(i) \in [N]$ denote the index of the medoid of the cluster containing $o_i$ at level $k$. The objective value is then
\begin{align*}
    \sum_{k\in [K]}  \sum_{i\in [N]} d_{i \mu_{k}(i)}.
\end{align*}

\section{Benders formulations}\label{section3:formulations}

In this section, we present Benders formulations based on the $p$-median literature. In the next section, we provide solution methods to solve the formulations and explain how to transform the obtained nested $p$-median solution into a hierarchical $k$-medoids solution.

\subsection{Standard Benders formulation}
In this section, we adapt a Benders formulation for the $p$-median problem from \cite{duranmateluna2023} to a nested $p$-median problem.

Let $\theta_{ik}$ be a decision variable denoting the distance from location $i$ at level $k$ to its nearest facility. Furthermore, let $y_{jk}$ be a binary variable indicating that facility $j$ is selected at level $k$. Let $N_i\leq N$ be the number of unique distances from $i$ to any location. We define the sorted distances $D_i^0\leq D_i^1 \leq \dots \leq  D_i^{N_i}$, where $D_i^0=-\infty$. A Benders formulation can be written as
\begin{align}
    \min & \sum_{k\in [K]} \sum_{i\in [N]} \theta_{ik}, \label{form2:obj}\\
    \text{s.t.} & \sum_{j\in [N]} y_{jk} = k, & k\in [K]\label{form2:cons1}\\
     & y_{jk}\leq y_{j(k+1)}, & j\in [N], k\in [K-1], \label{form2:cons4}\\
   & \theta_{ik} \geq D_i^{{n_i}+1} - \sum_{j: d_{ij}\leq D_i^{n_i}} (D_i^{{n_i}+1} - d_{ij}) y_{jk}, &   i\in [N], {n_i}\in [0, N_i-1], k\in [K],\label{form2:cons3}\\
   & y_{jk} \in \{0, 1\}, & j\in [N], k\in [K]. \label{form2:cons5}
\end{align}
The objective~\eqref{form2:obj} models the sum of $k$-medoids objectives across all levels. At each level $k$, exactly $k$ centers are selected, enforced by \eqref{form2:cons1}. The nestedness of the centers is enforced by \eqref{form2:cons4}. The Benders cuts are given by \eqref{form2:cons3}. The domain of the variables is enforced by \eqref{form2:cons5}.

\subsection{Aggregated Benders formulation}
One approach to addressing large-scale location problems is to perform aggregation. Aggregation can be applied to the problem data, for example by aggregating locations, or directly to the mathematical model through the aggregation of variables or constraints. Such aggregation techniques have been used for location problems in both heuristic approaches \citep{irawan2015a, avella2012, irawan2014, irawan2015b, salhi2015} and exact methods \citep{fischetti2016, ljubic2024, contreras2011, ramirez2023, willemsen2025accelerating}.

In this paper, we consider the aggregation of both constraints and variables, similar to \cite{willemsen2025accelerating}. Let $\bar{\theta}_i:=\sum_{k=1}^K \theta_{ik}$ be a decision variable representing the sum of the assignment distances for observation $i$ across all levels $[K]$. We introduce a binary decision variable $z_{jk}$, which equals $1$ when observation $j$ is first selected as a center at level $k$ and remains selected until level~$K$. This leads to the aggregated formulation
\begin{align}
    \min & \sum_{i\in [N]} \bar{\theta}_{i}, \label{form3:obj}\\
    \text{s.t.} & \sum_{j\in [N]} z_{jk} = 1, & k\in [K]\label{form3:cons1}\\
     & \sum_{k\in [K]} z_{jk} \leq 1, & j\in [N]\label{form3:cons4}\\
   & \bar{\theta}_{i} \geq K D_i^{{n_i}+1} - \sum_{j: d_{ij}\leq D_i^{n_i}} (D_i^{{n_i}+1} - d_{ij}) [\sum_{k\in [K]} (K-k + 1) z_{jk}], &   i\in [N], {n_i}\in [0, N_i-1],\label{form3:cons3}\\
   & z_{jk} \in \{0, 1\}, & j\in [N], k\in [K]. \label{form3:cons5}
\end{align}
The objective \eqref{form3:obj} minimizes the sum of distances over all levels implicitly. Constraints~\eqref{form3:cons1} ensure that at each level $k$ exactly one new center is selected. Each location $j$ is selected as a new center at most once, which is enforced by \eqref{form3:cons4}. The aggregated Benders cuts are given by \eqref{form3:cons3}. These inequalities are obtained by substituting $y_{jk}=\sum_{m=1}^k z_{jm}$ into the standard Benders cuts \eqref{form2:cons3} and summing them over all levels using the same distance threshold. Since a center introduced at level~$m$ remains selected until level $K$, the variable $z_{jm}$ appears $K-m+1$ times. The resulting formulation is a relaxation because it only includes aggregated cuts with a common distance threshold across levels, whereas the standard Benders formulation allows a different threshold at each level. Consequently, solving this relaxation does not generally guarantee an optimal nested $p$-median solution.

\section{Solution methods}\label{section4:solutionMethods}

The Benders formulations presented in the previous section are usually solved via Benders decomposition. We follow a procedure similar to that proposed by \cite{duranmateluna2023} for the single-level $p$-median problem.

\subsection{Main method}
Since there are too many (aggregated) Benders cuts, we consider only a subset of cuts and we also relax the integrality constraints of the binary variables, which results in a so-called restricted master problem (RMP). When solving the RMP, some of the Benders cuts are violated. By identifying the violated cuts and adding them to the RMP we exclude the current infeasible RMP solution. These violated cuts are identified using a subproblem (SP). We iteratively solve the RMP and SP  until no violated cuts are detected. For the standard formulation, exact separation ensures that this yields an optimal solution to the linear relaxation. We then apply branch-and-Benders-cut to obtain an optimal integer solution. For the aggregated formulation, we use heuristic separation and stop the procedure after solving the root node, obtaining a nested $p$-median solution through rounding without an optimality guarantee.

The main method to solve both Benders formulations is summarized in \Cref{alg:mainMethod}. For the aggregated formulation, the procedure terminates at line 11. The inputs are distances between all locations and also the unique sorted distances. The algorithm returns the best lower and upper bounds on the objective function. The main method keeps track of the best lower bounds ($LB$) and upper bounds ($UB$). To obtain feasible upper bounds we perform a rounding heuristic and we perform two speedups based on \cite{duranmateluna2023}, namely constraint reduction and reduced cost fixing.

\paragraph{Rounding heuristic}
We construct nested center sets by processing the levels in increasing order of $k$. At each level $k$, we retain all previously selected centers and add unselected locations with the largest, possibly fractional, values $y_{jk}$ until exactly $k$ centers are selected. Assigning each observation to the nearest selected center at each level yields a feasible nested $p$-median solution. For the aggregated formulation, we apply the same procedure using the fractional $z_{jk}$ values. 

\paragraph{Constraint reduction}
\cite{duranmateluna2023} noted that when obtaining a fractional solution, most generated Benders cuts are not saturated by the fractional solution. These cuts can be removed to reduce the problem size. For each location $i$ and level $k$, let $\hat{n}_{ik}$ denote the largest distance index corresponding to a saturated constraint. We remove all cuts for that location and level with distance indices greater than $\hat{n}_{ik}$.

\paragraph{Reduced cost fixing}
Following \cite{duranmateluna2023}, we perform reduced cost fixing after solving the root node to optimality, using its objective value as the lower bound $LB$ and the best feasible solution value as the upper bound $UB$. Let $rc_{jk}$ denote the reduced cost of variable~$y_{jk}$. We fix $y_{jk}$ to zero if $LB + rc_{jk} > UB$, and to one if $LB - rc_{jk} > UB$.

\begin{algorithm}[h]
\caption{Main method}
\label{alg:mainMethod}
\begin{algorithmic}[1]
\Require sorted distances $D_i^0 \leq D_i^1\leq D_i^2 \leq \dots \leq D_i^{N_i}$
\Ensure best lower and upper bounds $LB^*$ and $UB^*$
\While{true}
    \State solve the RMP to obtain a solution $\matr{y}$ and a lower bound $LB$
    \State $UB\leftarrow$ \texttt{RoundingHeuristic}($\matr{y}$)
    \State $UB^* \leftarrow \min \{UB^*, UB\}$ and $LB^*\leftarrow LB$
     \State solve the separation problem on $\matr{y}$ to obtain violated cuts
     \If{no violated cut is identified}
\State break
     \EndIf
    \State add the violated cuts to the RMP
\EndWhile
\State stop if the solution is integral
\State perform constraint reduction and reduced cost fixing
\State apply branch-and-Benders-cut 
\end{algorithmic}
\end{algorithm}
\FloatBarrier

The main method returns the objective value of a nested $p$-median solution. In order to obtain a solution to our hierarchical $k$-medoids problem, we perform a translation step. We do this at every step at which we compute an upper bound in case an earlier iteration found a better hierarchical $k$-medoids solution.

\paragraph{Obtaining a hierarchical $k$-medoids solution}
It is not obvious how to translate a nested $p$-median solution into a hierarchical $k$-medoids solution without solving an additional optimization problem. Therefore, we heuristically translate the nested $p$-median solution into a fully nested $k$-medoids solution, where both center sets and cluster assignments are nested. This means that any fully nested $k$-medoids solution is also a feasible solution to our hierarchical $k$-medoids problem.

Starting from the nested center sets, we construct the nested cluster assignments recursively. At the first level, each location is assigned to its closest center. At each subsequent level, each newly selected center is assigned to the parent cluster containing that center, and each location is then assigned to its closest center among those belonging to its parent cluster. This restriction ensures that clusters at each level are nested within clusters at the previous level. After this conversion, a retained center need not be a medoid of its resulting cluster. We retain these centers when evaluating the solution, so we obtain an upper bound on the objective value of the constructed hierarchy.

\subsection{Subproblem}
In the SP for the standard Benders formulation, we identify the Benders cuts using the same method as in \cite{duranmateluna2023}, but we solve the separation procedure for each level~$k$. 

For the aggregated Benders formulation, we could perform the following separation procedure to identify violated aggregated Benders cuts. We define $w_j:= \sum_{k\in [K]} (K-k + 1) z_{jk}$. Then, the aggregated Benders cut \eqref{form3:cons3} can be separated using a modification of the standard separation procedure. For each location~$i$ we solve a separation procedure where we identify the critical position $r$ satisfying
\begin{align*}
    \sum_{j:d_{ij}< D_i^r} w_{j}<K
\leq
\sum_{j:d_{ij}\leq  D_i^r} w_{j}.
\end{align*}
However, in our experiments, the resulting separation procedure and added cuts increase the solution time of the RMP. Therefore, we use a heuristic separation procedure for the aggregated Benders formulation. Define $u_j:= \sum_{k\in [K]}  z_{jk}$. For each location $i$ and for each value $\beta\in [K]$, we solve a separation procedure to identify a critical position $r_\beta$ satisfying
\begin{align*}
    \sum_{j:d_{ij}< D_i^{r_\beta}} u_{j}<\beta
\leq
\sum_{j:d_{ij}\leq  D_i^{r_\beta}} u_{j}.
\end{align*}
We examine these candidates in increasing order of $\beta$ and evaluate the aggregated Benders inequality at the obtained threshold $D_i^{r_\beta}$ and add the first violated inequality found. Note that this separation procedure is not guaranteed to identify violated aggregated Benders inequalities even when one exists.

\section{Computational experiments}\label{section5:results}

\subsection{Setup}
In all experiments, we use the Manhattan distance rounded to five decimal places, although the method can be applied to other distance matrices. Unless stated otherwise, we impose a time limit of 600 seconds. In the tables, we report the model gap as $\frac{UB - LB}{UB} \times 100$, where $LB$ and $UB$ are the lower and upper bounds on the nested $p$-median objective obtained during the Benders iterations for the corresponding formulation. However, we report the best hierarchical $k$-medoids objective value, which is the performance measure of interest. We also report the percentage difference between objective values obtained by the different methods.

We consider three datasets. The first set is obtained from \cite{thrun2020}, for which the true number of clusters is known for each instance. These instances are widely used for evaluating clustering methods. The second set consists of clustering instances preprocessed by \cite{hennig2022}, which we divide into small (at most 200 observations) and large instances (between 2,000 and 5,000 observations). The third set consists of traveling salesman problem instances, which we refer to as TSP instances \citep{reinelt1991}.

All computational experiments are performed on an AMD Rome 7H12 processor at 2.6 GHz with 16 cores, and, by default, 32 GB of RAM. To solve the mixed integer programming formulations we use the commercial solver CPLEX 20.1. We only use a single thread, as this provides better computational performance \citep{duranmateluna2023}.

\subsection{TU instances}
First, we consider instances from \cite{thrun2020}. We compare our standard Benders decomposition method with the column generation (CG) approach from \cite{willemsen2026exact}. The latter considers subsets of the TU instances of size $N=50$ and $N=100$, which are publicly available\footnote{\url{https://pubsonline.informs.org/doi/suppl/10.1287/ijoc.2024.0903}}. In our tables, we report the computation times as reported in \cite{willemsen2026exact}. The experimental setup differs in some aspects. In particular, they did not explicitly round the distances, imposed a time limit of 1 hour compared with our 10 minute limit, and conducted their experiments on a similar machine with 32 cores and 64 GB RAM, whereas we use 16 cores and 32 GB RAM.

The results in \Cref{tab:thrun50} are obtained on TU instances with $N=50$ observations. The CG approach from \cite{willemsen2026exact} requires on average 426 seconds, with one instance reaching the time limit of 1 hour. In contrast, our standard Benders decomposition solves all instances in an average of 1.2 seconds. Although the standard Benders decomposition solves the hierarchical $k$-medoids problem heuristically, it obtains solutions of high quality. The largest percentage difference from an optimal hierarchical $k$-medoids objective value is less than 5\%. 

\begin{table}[h]
\centering
\caption{Results on TU instances with $N=50$. We report the computation time in seconds, the best hierarchical $k$-medoids objective value, and the percentage difference between the objective values. The column generation times are reported as in \cite{willemsen2026exact}, using a time limit of 3600 seconds.}
\label{tab:thrun50}
\small
\begin{tabular}{lrrrrrrrrr}
\toprule
\multicolumn{2}{c}{instance} &  & \multicolumn{2}{c}{CG} &  & \multicolumn{2}{c}{standard Benders} &  & diff (in \%)     \\ \cline{1-2} \cline{4-5} \cline{7-8}
name                   & $K$   &  & obj        & time         &  & obj            & time          &  &          \\ \midrule
Hepta                  & 7   &  &  -          & 3600.0    &  & 427.2          & 8.3        &  & - \\
dietary\_survey\_IBS   & 2   &  & 3110.2     & 1.5       &  & 3110.2         & 0.5        &  & 0.00     \\
Tetra                  & 4   &  & 309.2      & 180.2     &  & 312.7          & 0.5        &  & 1.13     \\
Lsun3D                 & 4   &  & 320.2      & 338.2     &  & 323.2          & 0.5        &  & 0.94     \\
Leukemia               & 6   &  & 43858.5    & 85.6      &  & 44622.8        & 0.6        &  & 1.74     \\
TwoDiamonds            & 2   &  & 148.0      & 0.4       &  & 148.9          & 0.5        &  & 0.62     \\
Target                 & 6   &  & 259.2      & 9.2       &  & 269.6          & 0.6        &  & 4.03     \\
Atom                   & 2   &  & 176.0      & 32.2      &  & 176.0          & 0.5        &  & 0.00     \\
Chainlink              & 2   &  & 198.2      & 39.5      &  & 198.2          & 0.5        &  & 0.00     \\
WingNut                & 2   &  & 133.3      & 0.4       &  & 134.2          & 0.5        &  & 0.63     \\
EngyTime               & 2   &  & 142.3      & 398.6     &  & 144.5          & 0.5        &  & 1.54     \\ \midrule
avg                    &     &  &            & 426.0     &  &                & 1.2        &  &         \\ \bottomrule
\end{tabular}
\end{table}

Similar results are obtained for the TU instances with $N=100$ observations. The CG approach reaches the time limit on 6 out of 11 instances. In contrast, the standard Benders decomposition method solves all instances in less than 1 second on average. It is on average 3,000 times faster than the previously reported average runtime of CG. Moreover, among the instances for which the CG approach proves optimality, the maximum percentage difference between the objective values is less than 3\%. 

\begin{table}[h]
\centering
\caption{Results on TU instances with $N=100$. We report the computation time in seconds, the best hierarchical $k$-medoids objective value, and the percentage difference between the objective values. The column generation times are reported as in \cite{willemsen2026exact}, using a time limit of 3600 seconds.}
\label{tab:thrun100}
\small
\begin{tabular}{lrrrrrrrrr}
\toprule
\multicolumn{2}{c}{instance} &  & \multicolumn{2}{c}{CG} &  & \multicolumn{2}{c}{standard Benders} &  & diff (in \%)     \\ \cline{1-2} \cline{4-5} \cline{7-8}
name                   & $K$   &  & obj        & time         &  & obj            & time          &  &          \\ \midrule
Hepta                  & 7   &  &         -   & 3600.0    &  & 872.3          & 1.4        &  & - \\
dietary\_survey\_IBS   & 2   &  & 6180.0     & 1.3       &  & 6180.0         & 0.6        &  & 0.00      \\
Tetra                  & 4   &  &         -   & 3600.0    &  & 694.8          & 0.7        &  & - \\
Lsun3D                 & 4   &  &         -   & 3600.0    &  & 643.0          & 0.5        &  & - \\
Leukemia               & 6   &  & 90263.1    & 2307.2    &  & 91742.9        & 1.1        &  & 1.64      \\
TwoDiamonds            & 2   &  & 291.1      & 0.7       &  & 295.7          & 0.4        &  & 1.56      \\
Target                 & 6   &  &         -   & 3600.0    &  & 570.9          & 1.0        &  & - \\
Atom                   & 2   &  &        -    & 3600.0    &  & 364.0          & 0.4        &  & - \\
Chainlink              & 2   &  &        -    & 3600.0    &  & 403.2          & 0.4        &  & - \\
WingNut                & 2   &  & 279.1      & 0.7       &  & 285.1          & 0.5        &  & 2.15      \\
EngyTime               & 2   &  & 281.3      & 30.8      &  & 288.2          & 0.4        &  & 2.48      \\ \midrule
avg                    &     &  &            & 2176.4    &  &                & 0.7        &  &              \\ \bottomrule
\end{tabular}
\end{table}
\FloatBarrier

Finally, for completeness, \Cref{tab:thrunN} reports the results for the original TU instances using all observations. The standard Benders decomposition method solves these instances in an average of 82 seconds.

\begin{table}[h]
\centering
\caption{Results on original TU instances. We report the computation time in seconds and the best hierarchical $k$-medoids objective value.}
\label{tab:thrunN}
\small
\begin{tabular}{lrrrrr}
\toprule
\multicolumn{3}{c}{instance}            &  & \multicolumn{2}{c}{standard Benders} \\ \cline{1-3} \cline{5-6} 
name                 & $N$    & $K$ &  &     obj           & time          \\ \midrule
Hepta                & 212  & 7 &  & 1860.5         & 6.6        \\
dietary\_survey\_IBS & 400  & 2 &  & 24974.9        & 1.5        \\
Tetra                & 400  & 4 &  & 3077.1         & 62.9       \\
Lsun3D               & 404  & 4 &  & 2570.5         & 13.9       \\
Leukemia             & 554  & 6 &  & 530608.9       & 39.1       \\
TwoDiamonds          & 600  & 2 &  & 2365.4         & 19.4       \\
Target               & 770  & 6 &  & 4648.7         & 110.0      \\
Atom                 & 800  & 2 &  & 2981.7         & 19.6       \\
Chainlink            & 1000 & 2 &  & 4092.1         & 19.1       \\
WingNut              & 1000 & 2 &  & 3018.6         & 67.4       \\
EngyTime             & 4096 & 2 &  & 11804.8        & 539.0      \\ \midrule
avg                  &      &   &  &                & 81.7    \\ \bottomrule
\end{tabular}
\end{table}
\FloatBarrier

\subsection{HENNIG instances}
In the remaining experiments, we evaluate the scalability and performance of both the standard and aggregated Benders decomposition methods. 

We first consider a set of small instances from \cite{hennig2022}. \Cref{tab:hennig_KN} reports the results for computing the full hierarchies, corresponding to $K=N$. The standard Benders decomposition requires on average 291 seconds, compared with only 8 seconds for the aggregated Benders decomposition. However, the aggregated Benders decomposition obtains substantially worse objective values than the standard Benders decomposition, with an average difference of 47.5\%.

\begin{table}[h]
\centering
\caption{Results on small HENNIG instances with $K=N$. We report the computation time in seconds, the best hierarchical $k$-medoids objective value, the relative optimality gap for the nested $p$-median formulation, and the percentage difference between the objective values.}
\label{tab:hennig_KN}
\resizebox{\textwidth}{!}{%
\begin{tabular}{lrrrrrrrrrrrr}
\toprule
\multicolumn{3}{c}{instance}         &  & \multicolumn{3}{c}{standard Benders} &  & \multicolumn{3}{c}{aggregated Benders} &  & diff (in \%)  \\ \cline{1-3} \cline{5-7} \cline{9-11}
name             & $N$   & $K$   &  & obj        & time    & gap    &  & obj        & time  & gap    &  &       \\ \midrule
London       & 32  & 32  &  & 2433.9     & 0.8     & 0.00   &  & 3292.7     & 0.3   & 34.49  &  & 35.29 \\
USArrests    & 50  & 50  &  & 1842.0     & 13.6    & 0.00   &  & 2525.8     & 0.6   & 38.11  &  & 37.12 \\
StoneFlakes  & 79  & 79  &  & 10641.7    & 47.7    & 0.00   &  & 13648.7    & 1.1   & 38.57  &  & 28.26 \\
Augsburg     & 95  & 95  &  & 18587.5    & 34.9    & 0.00   &  & 26847.8    & 1.7   & 39.39  &  & 44.44 \\
BreastTissue & 106 & 106 &  & 30088.4    & 268.1   & 0.00   &  & 36319.7    & 2.4   & 27.96  &  & 20.71 \\
Alcohol      & 125 & 125 &  & 7922.7     & 142.7   & 0.00   &  & 14656.8    & 4.1   & 79.15  &  & 85.00 \\
Iris         & 150 & 150 &  & 7366.9     & 600.0   & 0.14   &  & 13356.7    & 7.0   & 63.87  &  & 81.31 \\
BankNotes    & 200 & 200 &  & 39643.5    & 600.0   & 0.42   &  & 57199.7    & 17.1  & 41.93  &  & 44.29 \\
Crabs        & 200 & 200 &  & 12932.8    & 600.0   & 1.13   &  & 21132.3    & 20.1  & 64.89  &  & 63.40 \\
Heart        & 200 & 200 &  & 83849.0    & 600.0   & 2.00   &  & 113627.8   & 25.5  & 43.56  &  & 35.51 \\ \midrule
avg              &     &     &  &            & 290.8   & 0.37   &  &            & 8.0   & 47.19  &  & 47.53 \\ \bottomrule
\end{tabular}
}
\end{table}

The results for the large HENNIG instances are reported in \Cref{tab:hennig_large}. We consider instances with large values of $N$ and small values of $K$. For these instances, the aggregated Benders decomposition remains faster than the standard Benders decomposition. Interestingly, the percentage difference in objective value is substantially smaller than for the small instances and the average becomes negative for $K=6$. These results suggest that the aggregated Benders decomposition performs particularly well on larger instances when the number of hierarchy levels is relatively small.

\begin{table}[h]
\centering
\caption{Results on large HENNIG instances with varying $K$. We report the computation time in seconds, the best hierarchical $k$-medoids objective value, the relative optimality gap for the nested $p$-median formulation, and the percentage difference between the objective values.}
\label{tab:hennig_large}
\small
\resizebox{\textwidth}{!}{%
\begin{tabular}{llrrrrrrrrrrr}
\toprule
\multicolumn{3}{c}{instance}    &  & \multicolumn{3}{c}{standard Benders} &  & \multicolumn{3}{c}{aggregated Benders}  &  & diff (in \%)   \\ \cline{1-3} \cline{5-7} \cline{9-11} 
$K$ & name                 & $N$    &  & obj      & time   & gap    &  & obj      & time   & gap   &  &    \\ \midrule
\multirow{6}{*}{3} & Letters          & 2000 &  & 67775.0    & 224.3  & 0.00    &  & 67646.8  & 9.1   & 3.31  &  & -0.19  \\
                   & Cardiotocography & 2126 &  & 98923.5    & 71.6   & 0.00    &  & 103588.5 & 14.8  & 6.49  &  & 4.72   \\
                   & ImageSeg         & 2310 &  & 58116.6    & 453.9  & 0.00    &  & 61272.2  & 12.6  & 11.18 &  & 5.43   \\
                   & Decathlon        & 2580 &  & 59172.9    & 600.0  & 0.71    &  & 59652.5  & 27.2  & 4.29  &  & 0.81   \\
                   & Bats             & 2677 &  & 371946.6   & 94.1   & 0.00    &  & 380893.6 & 29.5  & 5.13  &  & 2.41   \\
                   & Spam             & 4601 &  & 203536.7   & 408.2  & 0.00    &  & 204607.4 & 36.8  & 1.23  &  & 0.53   \\ \cline{2-13}
                   & avg                  &      &  &            & 308.7  & 0.12    &  &          & 21.7  & 5.27  &  & 2.28   \\ \bottomrule
\multirow{6}{*}{4} & Letters          & 2000 &  & 88135.0    & 277.8  & 0.00    &  & 88081.5  & 15.0  & 3.57  &  & -0.06  \\
                   & Cardiotocography & 2126 &  & 127963.1   & 600.0  & 0.24    &  & 133154.4 & 30.8  & 7.28  &  & 4.06   \\
                   & ImageSeg         & 2310 &  & 71365.9    & 274.2  & 0.00    &  & 72971.2  & 21.7  & 9.33  &  & 2.25   \\
                   & Decathlon        & 2580 &  & 77408.6    & 600.0  & 1.46    &  & 79485.3  & 49.3  & 7.07  &  & 2.68   \\
                   & Bats             & 2677 &  & 478510.3   & 363.8  & 0.00    &  & 485973.4 & 44.3  & 5.48  &  & 1.56   \\
                   & Spam             & 4601 &  & 268594.5   & 600.0  & 1.22    &  & 269897.5 & 55.0  & 1.54  &  & 0.49   \\ \cline{2-13}
                   & avg                  &      &  &            & 452.6  & 0.49    &  &          & 36.0  & 5.71  &  & 1.83   \\ \bottomrule
\multirow{6}{*}{5} & Letters          & 2000 &  & 107161.7   & 600.0  & 0.48    &  & 108723.1 & 24.2  & 4.45  &  & 1.46   \\
                   & Cardiotocography & 2126 &  & 157087.1   & 600.0  & 0.81    &  & 166083.5 & 39.5  & 8.45  &  & 5.73   \\
                   & ImageSeg         & 2310 &  & 82667.3    & 348.2  & 0.00    &  & 91022.2  & 26.1  & 15.96 &  & 10.11  \\
                   & Decathlon        & 2580 &  & 94450.4    & 600.0  & 1.81    &  & 99005.4  & 75.9  & 8.89  &  & 4.82   \\
                   & Bats             & 2677 &  & 579558.7   & 555.8  & 0.00    &  & 605932.0 & 55.9  & 9.00  &  & 4.55   \\
                   & Spam             & 4601 &  & 331386.0   & 600.0  & 0.80    &  & 343315.1 & 73.6  & 4.86  &  & 3.60   \\ \cline{2-13}
                   & avg                  &      &  &            & 550.7  & 0.65    &  &          & 49.2  & 8.60  &  & 5.04   \\ \bottomrule
\multirow{6}{*}{6} & Letters          & 2000 &  & 127336.3   & 600.0  & 1.19    &  & 131108.8 & 29.1  & 6.83  &  & 2.96   \\
                   & Cardiotocography & 2126 &  & 187785.2   & 600.0  & 1.90    &  & 196790.9 & 68.2  & 10.10 &  & 4.80   \\
                   & ImageSeg         & 2310 &  & 93621.8    & 543.1  & 0.00    &  & 102741.5 & 33.4  & 19.44 &  & 9.74   \\
                   & Decathlon        & 2580 &  & 110641.4   & 600.0  & 2.00    &  & 111852.1 & 136.5 & 5.82  &  & 1.09   \\
                   & Bats             & 2677 &  & 929437.3   & 600.0  & 100.00  &  & 711156.4 & 60.4  & 8.86  &  & -23.49 \\
                   & Spam             & 4601 &  & 399522.7   & 600.0  & 4.30    &  & 394783.9 & 92.6  & 1.62  &  & -1.19  \\ \cline{2-13}
                   & avg                  &      &  &            & 590.5  & 18.23   &  &          & 70.0  & 8.78  &  & -1.01  \\ \bottomrule
\end{tabular}
}
\end{table}
\FloatBarrier

\subsection{TSP instances}
The largest instances we consider are TSP instances. For these instances, we increase the available memory to 224 GB RAM and limit computation time to 30 minutes.

The results for $K=3$ are shown in \Cref{tab:TSP}. The standard Benders decomposition reaches the time limit for 2 out of 5 instances, leading to an average computation time of 1276 seconds. In contrast, the aggregated approach has an average computation time of 301 seconds. Furthermore, the aggregated Benders decomposition obtains objective values that are, on average, 9.5\% better than those obtained by the standard Benders decomposition. This average improvement is driven by the largest instance, for which the objective value is 65.8\% lower, whereas the other objective values are higher by at most 7\%.

\begin{table}[h]
\centering
\caption{Results on TSP instances with $K=3$. We report the computation time in seconds, the best hierarchical $k$-medoids objective value, the relative optimality gap for the nested $p$-median formulation, and the percentage difference between the objective values.}
\label{tab:TSP}
\resizebox{\textwidth}{!}{%
\begin{tabular}{lrrrrrrrrrrr}
\toprule
\multicolumn{2}{c}{instance} &  & \multicolumn{3}{c}{standard Benders} &  & \multicolumn{3}{c}{aggregated Benders} &  & diff (in \%)   \\ \cline{1-2} \cline{4-6} \cline{8-10}
name             & N         &  & obj       & time    & gap     &  & obj     & time   & gap   &  &        \\ \midrule
d2103        & 2103      &  & 8757.2    & 192.8   & 0.00    &  & 9342.4  & 23.2   & 14.31 &  & 6.68   \\
pcb3038      & 3038      &  & 13074.6   & 1400.8
  & 0.00    &  & 13596.9 & 53.6   & 14.81 &  & 4.00   \\
fl3795       & 3795      &  & 15355.5   & 1187.9  & 0.00    &  & 15898.2 & 72.4   & 10.12 &  & 3.53   \\
rl5934       & 5934      &  & 25161.7   & 1800.0  & 5.60    &  & 26203.0 & 213.9  & 12.77 &  & 4.14   \\
usa13509     & 13509     &  & 156342.0  & 1800.0  & 100.00  &  & 53526.6 & 1140.3 & 10.57 &  & -65.76 \\ \midrule
avg              &           &  &           & 1276.3
  & 21.12   &  &         & 300.7  & 12.51 &  & -9.48  \\ \bottomrule
\end{tabular}
}
\end{table}
\FloatBarrier

\section{Conclusion}\label{section6:conclusion}
We considered a hierarchical clustering problem formulated as a sum of $k$-medoids objectives. We showed that techniques from the $p$-median literature can be adapted to heuristically solve the hierarchical $k$-medoids problem.

We proposed two approaches: the first is based on a standard Benders decomposition, while the second uses aggregated Benders cuts. Our computational experiments demonstrate that the standard Benders decomposition can be three orders of magnitude faster than an exact column generation approach from the literature, while obtaining objective values within 5\% of an optimal solution. The aggregated Benders approach further improves computational performance on large instances and scales to instances with more than 13,000 observations.

Overall, our results demonstrate that, despite the differences between nested $p$-median and hierarchical $k$-medoids problems, algorithmic techniques developed for one problem can be effectively adapted to the other.

\renewcommand{\bibname}{References} 
\phantomsection
\addcontentsline{toc}{section}{References}
\begingroup
\bibliographystyle{apalike}
\bibliography{sample} 
\endgroup

\newpage
\phantomsection
\addcontentsline{toc}{section}{Appendix}


\setcounter{table}{0}
\renewcommand{\thetable}{A\arabic{table}}
\setcounter{figure}{0}
\renewcommand{\thefigure}{A\arabic{figure}}
\renewcommand{\thesubsection}{\Alph{subsection}} 

\end{document}